\documentclass[12pt]{amsart}
\usepackage[centertags]{amsmath}
\usepackage{amsfonts,amssymb,mathrsfs,float}
\usepackage{colonequals}
\usepackage{graphicx}
\usepackage{placeins}
\usepackage{enumerate}
\usepackage[latin1]{inputenc}
\usepackage{float}

  \newtheorem{theorem}{Theorem}
  \newtheorem{corollary}{Corollary}
  \newtheorem{proposition}{Proposition}
  \theoremstyle{remark}

  \newcommand{\ZZ}{\mathbb Z}
  
   \newcommand{\QQ}{\mathbb Q}

 \newcommand{\midd}{{\,|\,}}
 \newcommand{\rad}{{\rm rad}}

\begin{document}

\title{Cyclotomic primes}
\date{\today}

\author[Carl Pomerance]{Carl Pomerance}\let\thefootnote\relax\footnote{Mathematics Department, Dartmouth College,
Hanover, NH 03755, USA.  email: carlp@math.dartmouth.edu}

\address{Mathematics Department, Dartmouth College, Hanover, NH 03755, USA}
\email{carlp@math.dartmouth.edu}

\begin{abstract}
Mersenne primes and Fermat primes may be thought of as primes of the form $\Phi_m(2)$,
where $\Phi_m(x)$ is the $m$th cyclotomic polynomial.  This paper discusses the more
general problem of primes and composites of this form.
\end{abstract}

\maketitle

\section{Introduction}
\label{S:intro}

Studied since antiquity, we have the Mersenne primes.  These are prime numbers 1 less
than a power of 2, so of the form $2^n-1$.  To be prime it is necessary that $n=p$ is
prime, but this is not sufficient, e.g., $p=11$.  The first 4 of these primes were known
to Euclid and they played a key role in his work on perfect numbers.  We now know more than 50 
Mersenne primes, the largest 
at present being $2^p-1$ with $p=136{,}279{,}841$,  see \cite{G}.  Evidently it takes some
doing to check the primality of numbers this large!

It is widely believed that there are infinitely many Mersenne primes, and also 
infinitely many primes $p$ with $2^p-1$ composite.  Though both assertions are
still unsolved, there is a conditional proof of the second one based on the prime $k$-tuples
hypothesis:  If $p\equiv3\pmod 4$ is prime with $p>3$ and $q=2p+1$ is prime, then
$2^p-1$ is composite.  Indeed, the conditions imply that $q\equiv7\pmod8$ so that
$(2/q)=1$.  This implies that $q\midd2^{(q-1)/2}-1=2^p-1$, and the condition $p>3$
implies that $q<2^p-1$.  Thus $q$ is a proper divisor of $2^p-1$ implying the latter
is composite.  For example, $23\midd2^{11}-1$.  It remains to note that the prime $k$-tuples
hypothesis implies there are infinitely many primes $p\equiv3\pmod 4$ with $2p+1$ prime.

Also studied for centuries are the Fermat primes.  These are primes that are 1 more
than a power of 2, so of the form $2^n+1$.  To be prime (and $>2$) it is necessary that
$n$ itself is a power of 2.  Again, this is not sufficient.  Fermat knew
that $2^{2^k}+1$ is prime for $k=0,1,2,3,4$ and he conjectured that it is always
prime.  However, Euler showed that $641\midd2^{2^5}+1$.  It is now known that
$2^{2^k}+1$ is composite for all larger values of $k$ up to 32, and also some sporadic
larger values as well.  It is conjectured that all but finitely many are composite and that
perhaps $2^{2^4}+1$ is the largest Fermat prime.  Nothing has been proved here, even
conditionally.
 
What Mersenne primes and Fermat primes have in common is that they are 
{\it cyclotomic primes}.  These are primes of the form $\Phi_m(2)$, where $\Phi_m$
is the $m$th cyclotomic polynomial.  This is the minimal polynomial in $\QQ[x]$ for $e^{2\pi i/m}$
and it has degree $\varphi(m)$, Euler's function.  We have the twin identities:
\[
x^m-1=\prod_{d\,|\,m}\Phi_d(x),\quad \Phi_m(x)=\prod_{d\,|\,m}(x^d-1)^{\mu(m/d)}.
\]
Note that if $p$ is prime, then $\Phi_p(x)=(x^p-1)/(x-1)$, so that $\Phi_p(2)=2^p-1$.
Further $\Phi_{2^{k+1}}(x)=(x^{2^{k+1}}-1)/(x^{2^k}-1)=x^{2^k}+1$, so that
$\Phi_{2^{k+1}}(2)=2^{2^k}+1$.  We also have the Wagstaff primes, see \cite{O},
which are primes of the form $\Phi_{2p}(2)=(2^p+1)/3$, where $p$ is an odd prime.

So, a cyclotomic prime is a prime of the form $\Phi_m(2)$.  We can ask if there are
infinitely many of them and also if there are infinitely many numbers of this form that
are composite.  It turns out that there are infinitely many composites
for fairly trivial reasons.  The substance of this paper is to show that there are
infinitely many nontrivial composites.  We make this precise in the next section.

\section{Basics and statement of results}

Let
\[
\phi_m:=\Phi_m(2).
\]
We say a prime factor $p$ of $\phi_m$ is {\it primitive} if it does not divide any $\phi_k$ for
$k<m$.  Otherwise we say $p$ is {\it intrinsic}.  For an odd prime $p$ let $\ell(p)$ denote the
mutlplicative order of 2 in $(\ZZ/p\ZZ)^\times$.  We have that $p$ is a primitive prime
factor of $\phi_m$ if and only if $\ell(p)=m$.  Further, $\phi_m$ has an intrinsic prime factor
$p$ if and only if $m=p^j\ell(p)$ for some positive integer $j$, in which case $p$ is the largest
prime factor of $m$ and $p\,\|\,\phi_m$.
If $m$ is of this form, let $\delta_m=p$, and otherwise let $\delta_m=1$.  Thus, every prime
factor of
\[
\psi_m:=\phi_m/\delta_m
\]
is primitive.  

We know that for each $m\notin\{1,6\}$,
there is at least one primitive prime factor of $\phi_m$; this is Bang's theorem.  
The numbers $\psi_m$ are pairwise coprime and except for $m=1$ or 6, they
are all $>1$ (cf.\ \cite{PR}).

Due to the factorization
\[
4x^4+1=(2x^2+2x+1)(2x^2-2x+1),
\]
there is a further generic factorization of $\phi_m$ beyond $\delta_m\psi_m$
when $m\equiv 4\pmod 8$.  Note that 
\[
2^{4k+2}+1=4(2^k)^4+1=(2^{2k+1}+2^{k+1}+1)(2^{2k+1}-2^{k+1}+1),
\]
which leads to the factorization
\begin{align}
\label{eq:Au}
\phi_{8k+4}&=\gcd(\phi_{8k+4},2^{2k+1}+2^{k+1}+1)\gcd(\phi_{8k+4},2^{2k+1}-2^{k+1}+1)\\
&=:\phi_{8k+4}^+\phi_{8k+4}^-.\notag
\end{align}
By dividing out an intrinsic prime factor if it exists, we have
\[
\psi_{8k+4}=:\psi_{8k+4}^+\psi_{8k+4}^-.
\]
Further, this factorization is nontrivial for $k\ge3$, a result due to Schinzel \cite{S}.  The
factorization \eqref{eq:Au} is known under the name Aurifeuille, see \cite{B}.

Note that 
\begin{equation}
\label{eq:phiineq}
\phi_m\in[2^{\varphi(m)-1},2^{\varphi(m)+1}),
\end{equation}
 see \cite[Theorem 3.6]{H}, \cite[Theorem 4.3]{PR}.  Also, using \cite[eqs.~(13), (14)]{B}
 it is not hard to show that
 \begin{equation}
 \label{eq:Auineq}
 \phi_{8k+4}^+,\phi_{8k+4}^-\asymp 2^{\varphi(8k+4)/2},
 \end{equation}
 where the notation indicates the 2 items on the left side are of the same magnitude
 as the item on the right side.

To state our results, consider the sets
\[
C_1=\{\psi_m:m\not\equiv4\kern-5pt\pmod8\},\quad C_2=\{\psi_m:m\equiv4\kern-5pt\pmod8\}.
\]

\begin{theorem}\label{th:main}
For sufficiently large values of $x$, the
set $C_1$ contains more than $x^{3/5}$ composite numbers $\psi_m$ with $m\le x$
and the set $C_2$ contains
more than $x^{3/5}$  numbers $\psi_m$ with $m\le x$ that are not the product of two primes.
\end{theorem}

The exponent $3/5$ in the theorem is not optimal, this is discussed below.
The proof uses the deep result that for some $\theta$ with $1/2<\theta<1$,
there are infinitely many primes $p$ such that
$p-1$ has a large prime factor $q>p^\theta$.  We can also
prove a slightly stronger result conditional on the abc conjecture.
\begin{theorem}
\label{th:abc}
Assume the abc conjecture.  The set $C_1$ contains infinitely many numbers divisible
by at least $2$ distinct primes and the set $C_2$ contains infinitely many numbers divisible by
at least $3$ distinct primes.
\end{theorem}
We remark that the abc conjecture has been used for similar purposes in \cite{Si}
and \cite[Theorem 3]{BFGS}.  We also remark that this theorem gives an 
abc-conjecture-conditional solution of a problem of Schinzel \cite[p. 561]{S}.

Throughout the letters $p,q$ will always denote prime numbers.  We also let $P^+(n)$
denote the largest prime factor of $n>1$, and we let $P^+(1)=1$.

\section{An elementary approach}

Here we prove a somewhat weaker version of Theorem \ref{th:main} where the exponent
$3/5$ is replaced with $1/2$.  Let $x$ be large and consider primes $p\le x$.
 It follows from Erd\H os--Murty \cite[Theorem 1]{EM} that the number of
such primes $p$ with $\ell(p) >5 x^{1/2}(\log x)^2$ is $\sim x/\log x$.  For any positive
integer $d$ let $L_d(x)$ denote the set of primes $p\in(x/2,x]$ with 
\begin{align*}
&p\,\equiv\,3\kern-5pt\pmod4,\\
&\ell(p)=(p-1)/d,\\
&\ell(p)>5x^{1/2}(\log x)^2.
\end{align*}
It follows that
\begin{equation}
\label{eq:Lsum}
\sum_{d\le \frac15x^{1/2}/(\log x)^2}\#L_d(x)\sim \frac x{4\log x}.
\end{equation}
Thus, for large $x$ there is some number $d_0\le \frac15x^{1/2}/(\log x)^2$ with
\[
\#L_{d_0}(x) > x^{1/2}.
\]
Consider now the values of $m=\ell(p)=(p-1)/d_0$ for $p\in L_{d_0}(x)$.  They
are all distinct, bounded by $x$, and $\not\equiv4\pmod8$.  For large $x$,
$\psi_m$ is easily seen to be $>x$ (using $m>x^{1/2}$ and \eqref{eq:phiineq})
and $p\midd\psi_m$.  It follows that $\psi_m$ is composite.  Thus, $C_1$ contains
more than $x^{1/2}$ composite numbers $\psi_m$ with $m\le x$.

By changing 3 (mod 4) above to 5 (mod 8) and using \eqref{eq:Auineq}, the analogous
argument shows that at least one of $\psi_m^+,~\psi_m^-$ is composite, so
$C_2$ contains more than $x^{1/2}$ numbers $\psi_m$ with $m\le x$  which are not the
product of two primes.

\section{Proof of Theorem \ref{th:main}}

Denote by $\pi(x;d,a)$ the number of primes $p\le x$ with $p\equiv a\pmod d$.
Our principal tool is the following theorem.  Let $\theta=3/5$.
\begin{proposition}
\label{prop:BH}
We have
\[
\sum_{q>x^\theta}\pi(x;4q,2q+1)\log q\gg x\quad\hbox{and}\quad
\sum_{q>x^\theta}\pi(x;8q,4q+1)\log q\gg x.
\]
\end{proposition}
The analogous result for $\pi(x;q,1)$ is well known with varying values of ``$\theta$" in the
literature.  The current champions are Baker and Harman \cite{BH}, who essentially have
$\theta=0.677$, though they do not state their result in the same way.  Probably the techniques
of their paper would allow the same value of $\theta$ in Proposition \ref{prop:BH}, but we
do not pursue the optimal value at this point.  Other results in their paper have been
recently strengthened (see \cite{L}); conjecturally any value of $\theta <1$ may be used.

We now sketch a proof of Proposition \ref{prop:BH}.
With $\Lambda$ the von Mangoldt function, we have
\begin{align*}
\sum_{d\le x}\pi(x;4d,2d+1)\Lambda(d)
&=\sum_{\substack{p\,\equiv\,3\kern-5pt\pmod4\\p\le x}}\sum_{d\midd(p-1)/2}\Lambda(d)+O(x/\log x)\\
&=\sum_{\substack{p\,\equiv\,3\kern-5pt\pmod4\\p\le x}}\log(p-1)+O(x/\log x) \\
&=\frac12x+O(x/\log x).
\end{align*}
Further, by the Bombieri--Vinogradov theorem plus
a small additional argument using the Brun--Titchmarsh inequality
(see \cite{MV}) to clean up the boundary cases, we have
\[
\sum_{d\le x^{1/2}}\pi(x;4d,2d+1)\Lambda(d)\sim \frac14x,\quad x\to\infty.
\]
Thus,
\[
\sum_{d>x^{1/2}}\pi(x;4d,2d+1)\Lambda(d)\sim \frac14x,\quad x\to\infty.
\]
The contribution to this last sum when $d$ is composite is $o(x)$, so we have
\[
\sum_{q>x^{1/2}}\pi(x;4q,2q+1)\log q\sim\frac14x,\quad x\to\infty.
\]
By the Brun--Titchmarsh inequality,
 \begin{align*}
 \sum_{x^{1/2}<q\le x^{\theta}}\pi(x;4q,2q+1)\log q
 &\le\sum_{x^{1/2}<q\le x^{\theta}}\frac{2x\log q}{\varphi(4q)\log(x/4q)}\\
& \sim x\log(5/4)<0.23x.
 \end{align*}
Thus, with the prior display, we have
\[
\sum_{q>x^\theta}\pi(x;4q,2q+1)\log q \gg x,
\]
which shows the first assertion in Proposition \ref{prop:BH}.
The second assertion follows in a similar manner.

To achieve a cosmetically more appealing version of our results, note that since
\[
\sum_{x^\theta<q\le x^\theta(\log x)^2}\frac1q\ll\frac{\log\log x}{\log x}=o(1),
\]
the Brun--Titchmarsh theorem implies that in both parts of Proposition~\ref{prop:BH}
we may replace $q>x^\theta$ with $q>x^\theta(\log x)^2$.

Since no prime $p\le x$ is divisible by 2 different primes $q>x^\theta$, we have the 
following result.
\begin{corollary}
\label{cor:primecount}
We have
\[
\sum_{\substack{p\,\equiv\,3\kern-5pt\pmod4\\P^+(p-1)>x^\theta(\log x)^2\\p\le x}}1\gg x/\log x
\quad\hbox{and}\quad
\sum_{\substack{p\,\equiv\,5\kern-5pt\pmod8\\P^+(p-1)>x^\theta(\log x)^2\\p\le x}}1\gg x/\log x.
\]
\end{corollary}

An elementary argument shows that the
number of primes $p$ with $\ell(p)=k$ is $\ll k/\log k$, so it follows that the
number of primes $p\le x$ with 
$\ell(p)\le x^{1-\theta}/(\log x)^2$ is $\ll x^{2(1-\theta)}=o(\pi(x))$. (Note that there is
a similar argument in \cite{Gol}.  Also here one could appeal to
\cite[Theorem~1]{EM}.)  However, a prime $p$ counted in either part of
Corollary \ref{cor:primecount} either has $\ell(p)\le x^{1-\theta}/(\log x)^2$ or $\ell(p)>x^\theta(\log x)^2$.
Hence we have
\[
\sum_{\substack{p\,\equiv\,3\kern-5pt\pmod4\\\ell(p)>x^\theta(\log x)^2\\p\le x}}1\gg x/\log x
~\hbox{ and }~
\sum_{\substack{p\,\equiv\,5\kern-5pt\pmod8\\\ell(p)>x^\theta(\log x)^2\\p\le x}}1\gg x/\log x.
\]
Thus, using the notation of the previous section we have the improvement on
\eqref{eq:Lsum}:
\[
\sum_{d\le x^{1-\theta}/(\log x)^2}\#L_d(x)\gg\frac x{\log x},
\]
and an analogous result holds for primes $p\equiv5\pmod 8$.
Thus, by the same argument as in the previous section, we have Theorem \ref{th:main}.
   
 \section{Conditional results}
 
 For a positive integer $n$ let $\rad(n)$ denote the largest squarefree divisor of $n$.
 The abc conjecture asserts that for each fixed $\epsilon>0$, there are at most finitely
 many coprime positive integer triples $a,b,c$ with $a+b=c$ and $\rad(abc)<c^{1-\epsilon}$.
 In this section we will prove Theorem \ref{th:abc}, which is conditional on the abc conjecture,
  and also discuss some other conditional results.
 
 \begin{proof}[Proof of Theorem \ref{th:abc}]
 
 First suppose that $m\not\equiv4\pmod 8$.  We know from Theorem \ref{th:main} that
 there are infinitely many such $m$ with $\psi_m$ composite, and in fact, there are 
 $\ge x^\theta$ such $m\le x$ when $x$ is large.   Further, each such $m$ is of the form
 $k_qq$ where $q>x^\theta$ and $\psi_{k_qq}$ is divisible by a prime $p=p_q\le x$
with $\ell(p)=k_qq$.
 The only way for a composite number
 not to be divisible by at least 2 distinct primes is if it is a prime power, namely $p^i$ where
 $i\ge2$, so  suppose that $\psi_{k_qq}=p^i$.
 Consider the abc equation
 $1+(2^{k_qq}-1)=2^{k_qq}$.  Since $\psi_{k_qq}\ge\phi_{k_qq}/q$, we have
 \[
 \rad(abc)\le 2pq(2^{k_qq}-1)/\phi_{k_qq}.
 \]
 Assuming the abc conjecture this would be impossible for large $q$
 if there is some fixed $\epsilon>0$ such that $\phi_{k_qq}>2^{\epsilon k_qq}$
 (since $2pq=O(x^2)=2^{o(k_qq)}$).
 Using \eqref{eq:phiineq} this would follow if $\varphi(k_q)>2\epsilon k_q$.  We now
 show that we may assume this is indeed the case.
 
 Recalling that $\ell(p)=k_qq$, write $j=(p-1)/q$, so that $k_q\midd j$.  And so if
 $\varphi(k_q)/k_q\le2\epsilon$, then $\varphi(j)/j\le 2\epsilon$.  For a given value of $j\le x^{1-\theta}$
 we consider primes $q\le x/j$ such that $jq+1$ is
 prime.  By Brun's or Selberg's sieve, the number of such $q$ is $\ll x/(\varphi(j)(\log x)^2)$.
 In the proof of Proposition \ref{prop:BH} we showed there are $\ge(.02+o(1))x/\log x$
 pairs $q,p$.  Let
 \[
 J=\{j\le x^{1-\theta}:\varphi(j)/j\le2\epsilon\}.
 \]
 We will show that
$\sum_{j\in J}1/\varphi(j)\ll\epsilon\log x$ with $\epsilon$ small enough this will be
negligible in comparison with $(.02+o(1))\log x$.  This would follow from Erd\H os \cite[Theorem 1]{E}
(also see \cite[Theorem B]{KPP}), but we prefer to use the simpler approach in
\cite[Section 3]{KKP}.
 
  We have $(j/\varphi(j))^2=\sum_{d\midd j}h(d)$, where $h$ is multiplicative, supported
  on the squarefrees, and has $h(p)=(2p-1)/(p-1)^2$.  Then
  \[
  \sum_{n\le z}\Big(\frac n{\varphi(n)}\Big)^2=\sum_{d\le z}h(d)\left\lfloor\frac zd\right\rfloor
  <z\prod_{p}\Big(1+\frac{h(p)}p\Big)<4.5z.
  \]
  Thus, for any $\delta>0$,
  \[
  \sum_{\substack{n\le z\\\varphi(n)/n\le\delta}}1<4.5\delta^2z.
  \]
  A partial summation argument then shows that
  \[
  \sum_{\substack{n\le z\\\varphi(n)/n\le\delta}}\frac1n<4.5\delta^2\log(ez),
  \]
  so that writing $1/\varphi(n)=1/n\cdot n/\varphi(n)$,
  \[
  \sum_{\substack{n\le z\\\frac12\delta<\varphi(n)/n\le\delta}}\frac1{\varphi(n)}
  <\frac2\delta4.5\delta^2\log(ez)=9\delta\log(ez).
  \]
  We apply this with $z=x^{1-\theta}$ and at $\delta=2\epsilon,\epsilon,\frac12\epsilon,\dots$
  getting
  \[
  \sum_{j\in J}\frac1{\varphi(j)}<36\epsilon\log(ex).
  \]
  Thus if $\epsilon=\epsilon_0$ is small enough, we would have the number of $q,p$
  pairs with $\varphi(j)/j \le\epsilon_0$ being $<.01x/\log x$.  But we have seen in the
  previous section that the total number of $q,p$ pairs generated is $\ge(.02+o(1))x/\log x$.
  Thus, we can discard those with $\varphi(j)/j \le\epsilon_0$ and still be left with
  $\gg x/\log x$ pairs.  So, we will have $\varphi(j)/j>\epsilon_0$, which implies that
  $\varphi(k_q)/k_q>\epsilon_0$.
Thus, the abc conjecture is in play to show that the equation $1+(2^{k_qq}-1)=2^{k_qq}$
 with $\psi_{k_qq}$ a power of $p$ cannot occur when $x$ is large.
  
  The situation for $kq\equiv4\pmod 8$ is completely analogous; we suppress the details.
  \end{proof}
   
   By being more careful with the estimates one can get counts in Theorem \ref{th:abc}
  like the ones we found in Theorem \ref{th:main}.
  We further remark that a variant of our proof can show that for asymptotically all $m$,
  $\psi_m$ is not square-full, and the same goes for $\psi_{8k+4}^+$ and $\psi_{8k+4}^-$.
  
  We mentioned in the introduction that the prime $k$-tuples conjecture can be used
  to show that there are infinitely many primes $p$ with $\psi_p=2^p-1$ composite.
  We add here a couple of thoughts.  First, since we know that  8 and 9 form the only pair of consecutive
  numbers which are nontrivial powers (a result of Mih\u{a}ilescu \cite{M}), 
  it follows that $2^p-1$ cannot be a nontrivial power, so in this case, the abc conjecture
  is not necessary.  Second, using the Hardy--Littlewood version of the $k$-tuples conjecture,
  we have the number of primes $p\le x$ with $2^p-1$ divisible by at least 2 different primes
   is $\gg x/(\log x)^2$.
  
  We can prove there are more cyclotomic composites
assuming Artin's primitive root conjecture.  If
  2 is a primitive root for $p$, we have $p$ a prime factor of $\psi_{p-1}$ and so
  $\psi_{p-1}$ is composite for $p$ large.  By
  Hooley's GRH conditional proof of Artin's conjecture, we have $\gg x/\log x$
  primes $p\le x$ which have 2 as a primitive root.  Further, the proof is amenable
  to insisting that $p\equiv3\pmod 4$ and also the same holds when $p\equiv5\pmod 8$.
  So the GRH implies there are quite a few cyclotomic composites.  And as above, the
  abc conjecture can be used to show these composites are usually not prime powers.
  This result can be improved a little by considering primes $p \le kx$ with $\ell(p)=(p-1)/k$
  for various small values of $k$ and using sieve methods to show that $(p-1)/k=(p'-1)/k'$
  has few solutions when $k\ne k'$ are small.  Thus, with a little work it may be possible
  to show, assuming the GRH, that there are $\gg x\log\log x/\log x$ integers $l\le x$
  of the form $\ell(p)$ for some prime $p\ll x\log x$.
  
\section{Statistics and surmises}

  Concerning Table \ref{Ta:counts},  Gallot \cite{Gal} previously enumerated the cases
where $\phi_m$ is prime for $m\le6500$ and Noe \cite{N} extended this to $10^5$.  Our calculations agree with theirs.  In our work we used
 Mathematica and in particular their PrimeQ function.  This function is discussed in
 \cite{BFW}, where it is said to be based on the Baillie--PSW primality test.
 This is not a rigorous primality test, though no counterexamples are known (and there
 is a reward for the first one to be identified).  In fact, I have a heuristic argument that
 there are indeed infinitely many counterexamples, see \cite{Pdopo}.
So, it is possible that some of the prime
declarations made are false, but this seems unlikely, given that there are not very many
of them.   One of the larger primes unearthed here is $\phi_{60{,}287}$ which has 
$17{,}090$ decimal digits.  Note that when PrimeQ declares a number is not prime, this
conclusion is not in doubt.  Since PrimeQ is notably slower than checking if the
Fermat congruence $3^n\equiv3\pmod n$ holds, we first used that and confirmed
the few primality assertions with PrimeQ.  (We used the base 3 since every $\psi_m$ 
is either a prime or a base 2 pseudoprime.  See \cite{Po} where these thoughts are
developed.)  Many of the large primes uncovered here have indeed been certified (including
$\phi_{60{,}287}$) in the ongoing project factordb.com.  (Thanks are due to Yves Gallot
for informing me about this.)

  \begin{footnotesize}
\begin{table}[h!]
\caption{Counts for $m\le2^k$ with $\phi_m$ prime,  $\psi_m$ prime, $\psi_m^+$ prime, $\psi_m^-$ prime}
\label{Ta:counts}
\begin{tabular}{|rrcccc|} \hline
$k$&$\#m$ with $\phi_m$ prime & $\psi_m$ prime&$\psi_m^+$ prime&$\psi_m^-$ prime&\\ \hline
1	&	$1\kern20pt$&1&0&0&\\
2	&	$3\kern20pt$&3&1&0&\\
3	&	$7\kern20pt$&6&1&0&\\
4	&	$14\kern20pt$&13&2&0&\\
5	&	$23\kern20pt$&25&4&1&\\
6	&	$33\kern20pt$&36&7&5&\\
7	&	$49\kern20pt$&52&13&8&\\
8	&	$64\kern20pt$&68&20&16&\\
9	&	$81\kern20pt$&86&24&25&\\
10	&	$99\kern20pt$&106&30&33&\\
11	&	$122\kern20pt$&129&34&43&\\
12	&	$140\kern20pt$&147&44&54&\\
13	&	$167\kern20pt$&174&50&59&\\
14	&	$195\kern20pt$&202&61&64&\\
15	&	$221\kern20pt$&228&72&74&\\
16    &     $ 255\kern20pt$&262&85&83&\\
17    &     $289\kern20pt$&296&96&94&\\
\hline
\end{tabular}
\end{table}
\end{footnotesize}

Heuristically there are at most finitely many examples where $\psi_{p^i\ell(p)}$, with
$i\ge1$, is prime.
Is $\psi_{127\cdot7}$ the largest such example?  It is a prime of 226 decimal digits.
There are several examples where
both $\psi_m^+$ and $\psi_m^-$ are prime.  The largest that we found in our calculations to
$2^{17}$ is $m=1132$, where the two primes each have 85 decimal digits.
Probably there are at most finitely many of these ``twin cyclotomic primes".  

The counts
in Table \ref{Ta:counts} look to be proportional to $k^2$, and this is supported heuristically
as well.  Indeed, one can model $\psi_m$ as a random number near $2^{\varphi(m)}$
which has all prime factors larger than $m$.  So the ``probability" that it is prime
(given that $m\not\equiv4\pmod8$) is about $e^\gamma\log m/\varphi(m)\log2$.
The sum of these quantities up to $2^{15}$ is about 223.4,  up to $2^{16}$ is about
254.4, and up to $2^{17}$ is about 287.4, which are not bad matches with the table.
It would seem that the counts in the first column are asymptotically equal to $k^2$,
but this is likely not true.
One can sum $e^\gamma\log m/\varphi(m)\log2$ for $m\le 2^k$ with $m\not\equiv4\pmod8$,
finding it to be $\sim ck^2$, where 
\[
c=\frac5{12}e^\gamma\zeta(2)\zeta(3)\zeta(6)^{-1}\log 2=0.999774\dots.
\]
So, close to 1, but not 1.

One can enlarge further the realm of cyclotomic primes to look at the primitive parts of
$a^n-1$, where $a>2$.  Also one can look at the Fibonacci sequence, as well as other
Lucas sequences, for example see Drobot \cite{D}.
We suspect our methods carry over, but we leave this topic for another day, and perhaps another person.   
 
 \section*{Acknowledgments}
 I thank Max Alexeev, Michael Filaseta, Yves Gallot, Florian Luca,
 Mits Kobayashi, Pieter Moree,
  Paul Pollack, and Sam Wagstaff for their valuable comments and interest.

\section*{Appendix}

Here we list the values of $m$ corresponding to the counts in Table \ref{Ta:counts}.
\smallskip

\noindent
{\bf Values of $m$ with $\phi_m$ prime:}
\medskip

\noindent
2,\ 3,\ 4,\ 5,\ 6,\ 7,\ 8,\ 9,\ 10,\ 12,\ 13,\ 14,\ 15,\ 16,\ 17,\ 19,\ 22,\ 24,\ 26,\ 27,\ 30,\ 31,\ 32,\ 33,\  
 34,\  38,\  40,\  42,\  46,\  49,\  56,\  61,\  62,\  65,\  69,\  77,\  78,\  80,\  85,\  86,\  
 89,\  90,\  93,\  98,\  107,\  
120,\  122,\  126,\  127,\ 
129,\  133,\  145,\  150,\  158,\  165,\  170,\  174,\  184,\  192,\  195,\  202,\  
208,\   234,\  
254,\  261,\  280,\  296,\  312,\  322,\  334,\  345,\  366,\  374,\  382,\  398,\  410,\  414,\  
425,\  447,\  471,\  507,\  521,\  550,\  567,\  579,\  590,\  600,\  607,\  626,\  690,\  694,\  
712,\  745,\  
795,\  816,\  897,\  909,\  954,\  990,\ 
1106,\  1192,\  1224,\  1230,\  1279,\  1384,\  1386,\  1402,\ 
 1464,\  1512,\  1554,\  
1562,\  1600,\  1670,\  1683,\  1727,\  1781,\  1834,\  1904,\  1990,\  1992,\ 
2008,\  
2037,\  2203,\  2281,\  2298,\  2353,\  2406,\  2456,\  2499,\  2536,\  2838,\  3006,\  
3074,\  
3217,\  3415,\  3418,\  3481,\  3766,\  3817,\  3927,\ 4167,\ 4253,\ 4423,\ 4480,\ 5053,\ 5064,\ 5217,\ 5234,\ 5238,\ 5250,\ 5325,\ 
5382,\ 5403,\ 5421,\ 6120,\ 6925,\ 7078,\ 7254,\ 7503,\ 7539,\ 7592,\ 7617,\ 
7648,\ 7802,\ 7888,\ 7918,\ 8033,\ 
8370,\  9583,\  9689,\  9822,\  9941,\  10192,\ 
10967,\  11080,\  11213,\  11226,\  
11581,\  11614,\  11682,\  11742,\  11766,\  12231,\  12365,\  
12450,\  12561,\  13045,\  
13489,\  14166,\  14263,\  14952,\  14971,\  15400,\  15782,\  15998,\ 
16941,\  17088,\  17917,\  18046,\  19600,\  19937,\  20214,\  20678,\  21002,\  21382,\  
21701,\ 
 22245,\  22327,\  22558,\  23209,\  23318,\  23605,\  23770,\  24222,\  24782,\  
27797,\  28958,\ 
 28973,\  29256,\  31656,\  31923,\ 
 33816,\  34585,\  35565,\  35737,\  36960,\  39710,\  40411,\ 
  40520,\  42679,\  42991,\  
43830,\  43848,\  44497,\  45882,\  46203,\  47435,\  48387,\  48617,\  
49312,\  49962,\  
49986,\  50414,\  51603,\  51945,\  53977,\  55495,\  56166,\  56898,\  56955,\  
57177,\  
58315,\  58534,\  58882,\  60287,\ 67235,\ 67854,\ 69933,\ 70129,\ 70617,\ 75302,\ 76912,\ 78077,\ 
78426,\ 80160,\ 81165,\ 81432,\ 82569,\ 82730,\  84897,\ 85474,\ 85881,\ 86243,\ 87005,\ 94914,\ 
 95349,\ 99992,\ 100917,\ 104550,\ 108535,\ 109965,\ 110503,\ 110845,\ 111065,\ 116629,\ 
118080,\ 119210,\ 121806,\ 130002

\bigskip

\noindent
{\bf Values of $m$ with $\psi_m<\phi_m$ and $\psi_m$ prime:}
\medskip

\noindent
18,\ 20,\ 21,\ 54,\ 147,\ 342,\ 602,\ 889

 \bigskip
 
 \noindent 
{\bf Values of $m$ with $\psi_m^+$ prime:}
\medskip

\noindent
4,\ 12,\ 20,\ 28,\ 
36,\  44,\  60,\  68,\  76,\  84,\  100,\  108,\  116,\  132,\  140,\  180,\  204,\  220,\  228,\  252,\  
276,\  340,\  356,\  484,\  588,\  628,\  652,\  700,\  756,\  924,\  1132,\  1292,\  
1452,\  1516,\  2300,\  2484,\  2604,\  2964,\  3116,\  3276,\  3420,\  3540,\  
3940,\  3988,\  4892,\  5100,\  5268,\  5908,\  6620,\  7812,\  8964,\  9084,\  9324,\  
9468,\  10308,\  11980,\  12188,\  12204,\  13724,\  13860,\  15252,\  17052,\  18476,\  
20676,\  21916,\  24252,\  25004,\  25508,\  28692,\  29460,\  29492,\  31692,\  34236,\  
34380,\  35700,\  38428,\  40564,\  41316,\  45028,\  46076,\  50332,\  51148,\  51204,\  
56588,\  58796,\ 73668,\ 81900,\ 84020,\ 86508,\ 87420,\ 92324,\ 96204,\ 97524,\ 97620,\ 
104620,\ 118748
 
 \bigskip

\newpage
\noindent 
{\bf Values of $m$ with $\psi_m^-$ prime:}
\medskip

\noindent
28,\  36,\  44,\  52,\  60,\  84,\  108,\  116,\  132,\  140,\  172,\  188,\  196,\  212,\  220,\  
252,\  260,\  276,\  292,\  316,\  348,\  372,\  420,\  444,\  452,\  516,\  604,\  668,\  812,\  868,\  
924,\  956,\  964,\  1044,\  1132,\  1204,\  1276,\   1412,\  1468,\  1500,\  
1540,\  1564,\  1828,\  2124,\  2172,\  2228,\  2252,\  2452,\  2532,\  2716,\  2764,\  
2868,\  3484,\  3852,\  4844,\  5316,\  5468,\  6164,\  7828,\  9516,\  9684,\  10924,\  
12164,\  15860,\  19516,\  20588,\  21292,\  24180,\  25100,\  25212,\  28612,\  30988,\  
31460,\  32340,\  34404,\  38132,\  42660,\  43084,\  46292,\  46980,\  52740,\  56668,\  
60676,\ 68748,\ 69828,\ 72948,\ 74220,\ 75484,\ 84900,\ 87940,\ 106412,\ 110116,\ 
115292,\ 129836

\bigskip
\bigskip
\end{document}